\documentclass{article}
\usepackage{mathrsfs}
\usepackage{amsmath}
\usepackage{amssymb, amsmath}
\usepackage{graphicx}
\catcode`\@=11 \@addtoreset{equation}{section}

\catcode`\@=12
\usepackage{color}

\newtheorem{Theorem}{Theorem}[section]
\newtheorem{Proposition}{Proposition}[section]
\newtheorem{Lemma}{Lemma}[section]
\newtheorem{Corollary}{Corollary}[section]
\newtheorem{Definition}{Definition}[section]

\newtheorem{Remark}{Remark}[section]

\newcommand{\newcom}{\newcommand}
\newcommand{\bTheorem}[1]{
\begin{Theorem} \label{T#1} }
\newcommand{\eT}{\end{Theorem}}

\newcommand{\bProposition}[1]{
\begin{Proposition} \label{P#1}}
\newcommand{\eP}{\end{Proposition}}

\newcommand{\bLemma}[1]{
\begin{Lemma} \label{L#1} }
\newcommand{\eL}{\end{Lemma}}

\newcommand{\bCorollary}[1]{
\begin{Corollary} \label{C#1} }
\newcommand{\eC}{\end{Corollary}}

\newcommand{\beq}{\begin{equation}}
\newcommand{\eeq}{\end{equation}}
\newcom{\ben}{\begin{eqnarray}}
\newcom{\een}{\end{eqnarray}}
\newcom{\beno}{\begin{eqnarray*}}
\newcom{\eeno}{\end{eqnarray*}}
\newcom{\bali}{\begin{aligned}}
\newcom{\eali}{\end{aligned}}

\newcommand{\bFormula}[1]{
\begin{equation} \label{#1}}
\newcommand{\eF}{\end{equation}}

\newcommand{\f}{\frac}

\newcommand{\p}{\partial}

\newcommand{\vr}{\varrho}

\newcommand{\dy}{{\rm d} y}

\newcommand{\ds}{{\rm d} s}

\newcommand{\ep}{\varepsilon}

\newcommand\Cbox[2]{%
    \newbox\contentbox%
    \newbox\bkgdbox%
    \setbox\contentbox\hbox to \hsize{%
        \vtop{
            \kern\columnsep
            \hbox to \hsize{%
                \kern\columnsep%
                \advance\hsize by -2\columnsep%
                \setlength{\textwidth}{\hsize}%
                \vbox{
                    \parskip=\baselineskip
                    \parindent=0bp
                    #2
                }%
                \kern\columnsep%
            }%
            \kern\columnsep%
        }%
    }%
    \setbox\bkgdbox\vbox{
        \color{#1}
        \hrule width  \wd\contentbox %
               height \ht\contentbox %
               depth  \dp\contentbox
        \color{black}
    }%
    \wd\bkgdbox=0bp%
    \vbox{\hbox to \hsize{\box\bkgdbox\box\contentbox}}%
    \vskip\baselineskip%
}
\newcommand{\eq}[1]{\begin{equation}
\begin{split}
#1
\end{split}
\end{equation}}
\newcommand{\eqh}[1]{\begin{equation*}
\begin{split}
#1
\end{split}
\end{equation*}}

\begin{document}

%%%%%%%%%%%%%%%%%%%%%%%%%%%%%%%%

\title{\bf Large time behavior for a compressible two-fluid model with algebraic pressure closure and large initial data }

\author{Yang Li \\ Department of Mathematics, \\ Nanjing University, Nanjing 210093, China \\ lymath@smail.nju.edu.cn\\
Yongzhong Sun \\ Department of Mathematics, \\ Nanjing University, Nanjing 210093, China \\ sunyz@nju.edu.cn\\
Ewelina Zatorska \\ Department of Mathematics, \\ University College London, Gower Street, London WC1E 6BT,  United Kingdom \\
e.zatorska@ucl.ac.uk}

\maketitle
{\centerline {\bf Abstract }}
{In this paper, we consider a compressible two-fluid system with a common velocity field and algebraic pressure closure in dimension one. Existence, uniqueness and stability of global weak solutions to this system are obtained with arbitrarily large initial data. Making use of the uniform-in-time bounds for the densities from above and below, exponential decay of weak solution to the unique steady state is obtained without any smallness restriction to the size of the initial data. In particular, our results show that degeneration to single-fluid motion will not occur as long as in the initial distribution both components are present at every point. }



\section{Introduction}
In multi-component flows the presence of topologically  complex interphase separating the components is a great difficulty from physical as well as mathematical point of view. However, in most of engineering applications precise description of motion of each component of interphase are not rarely  needed and only the averaged macroscopic description is important. We will focus on the averaged two-component model derived in the monograph of Ishii and Hibiki in its inviscid form \cite{IsHi}. We refer the interested reader to \cite{BDGGH} for concise overview of various modelling and mathematical aspects related to such models. In the present paper we immediately assume that the two components of the flow share a common velocity field and that their pressures are equal (algebraic pressure closure). We obtain the following system of partial differential equations in dimension one:
\beq\label{bn1}
\p_t (\alpha_{\pm} \vr_{\pm})+ \p_x(\alpha_{\pm} \vr_{\pm} u)=0,
\eeq
\beq\label{bn2}
\p_t ( (\alpha_{+} \vr_{+}+\alpha_{-}\vr_{-})u   )+ \p_x((\alpha_{+} \vr_{+}+\alpha_{-}\vr_{-})u^2)+\p_x p=\mu \p_{xx} u,
\eeq
\beq\label{bn3}
\alpha_{+}+\alpha_{-}=1,\,\,\,\alpha_{\pm}\geq 0,
\eeq
\beq\label{bn4}
p=p_{+}=p_{-},
\eeq
where $\alpha_{+}$ and $\alpha_{-}$ are the volumetric rates of the two fluids; $\vr_{+}$ and $\vr_{-}$ are the two mass densities; $u$ is  the common velocity field, and $\mu>0$ is the viscosity coefficient. The two internal pressures are given by
\beq\label{bn5}
p_{+}=\vr_{+}^{\gamma_{+}  },\,\,\,p_{-}=\vr_{-}^{\gamma_{-}  },
\eeq
for adiabatic exponents $\gamma_{\pm}>1$. Following \cite{BMZ}, we introduce the notation
\beq\label{notRQ}
R=\alpha_{+} \vr_{+},\,\,\,Q=\alpha_{-}\vr_{-},\,\,\,Z=\vr_{+},
\eeq
and reformulate (\ref{bn1})-(\ref{bn5}) to
\beq\label{bn6}
\p_t R+\p_x (R u)=0,
\eeq
\beq\label{bn7}
\p_t Q+\p_x (Q u)=0,
\eeq
\beq\label{bn8}
\p_t ( (R+Q)u ) +\p_x ( (R+Q)u^2 ) +\p_x Z^{ \gamma_{+}  }=\mu \p_{xx}u.
\eeq
Due to the algebraic closure \eqref{bn4}, $Z$ is an implicit function of $R$ and $Q$ interrelated by
\beq\label{bn9}
Q=\left(1-\f{R}{Z}\right) Z^{\gamma},\,\,\, \gamma:= \f{\gamma_{+}}{\gamma_{-}},
\eeq
\beq\label{bn10}
R \leq Z.
\eeq
The same model, but in semi-stationary Stokes regime, has been recently investigated by Bresch, Mucha and the third author in the three-dimensional setting. They proved the global-in-time existence of weak solutions without any restriction on the initial data.
Similar result for the general Navier-Stokes system, with generalized equation of state was later obtained  by Novotn\'y and Pokorn\'y \cite{NP}. Earlier results in this spirit concern existence of weak solutions to very particular two-component models including the fluid model of atmospheric flow with transport of potential temperature \cite{MMM}, and the hydrodynamic limit of Vlasov-Fokker-Planck  system modelling suspension of the particles in the compressible fluid \cite{VWY}. For  other related results in case of one-dimensional two-fluid models, including density dependent coefficients or the so-called drift-flux model, we refer to \cite{BrLi}, to the works of Evje et. al \cite{EW1,EW2,EW3,EW4,EWZ1} and to the recent overview paper \cite{WYZ}.

The present paper is, as far as we know, the first attempt to provide some more information about quantitative  properties of  weak solutions to this system.
In order to investigate the large time behavior of solutions, we furthermore rewrite (\ref{bn6})-(\ref{bn8}) in Lagrangian coordinates. To do this, we assume that, for definiteness, the fluids occupy the closed interval $[0,1]$ with no-slip boundary conditions and make the change of variables
\[
y:=\int_0^x (R+Q)(\xi,t)d\xi,\,\,s:=t.
\]
As a consequence,
\beq\label{bn11}
\p_t \tau= \p_y u,
\eeq
%\beq\label{bn12a}
%\p_t (R \tau)=0,
%\eeq
\beq\label{bn12}
\p_t (Q \tau)=0,
\eeq
\beq\label{bn13}
\p_t u = \p_y \left(\mu \tau^{-1} \p_y u -Z(R,Q)^{ \gamma_{+}  } \right),
\eeq
where we set
\[
%\vr:= R+Q,\,\,\,
\tau:= (R+Q)^{-1}.
\]
The equations (\ref{bn11})-(\ref{bn13}) are supplemented with the initial and boundary conditions as follows:
\beq\label{bn14}
(R,Q,u)(y,0)=(R_0,Q_0,u_0)(y),\,\,y\in [0,1],
\eeq
\beq\label{bn15}
u(0,t)=u(1,t)=0,\,\,t\in [0,\infty),
\eeq
and we denote
\beq\label{bn14a}
\tau_0:=(R_0+Q_0)^{-1}.
\eeq

This paper is mainly devoted to the large time behavior of weak solutions to (\ref{bn11})-(\ref{bn15}) with large initial data. Existence, uniqueness and stability of weak solutions are obtained by making full use of the specific structure of the equations. Unlike in the three-dimensional regime \cite{BMZ, NP}, we prove the existence of weak solutions by approximation based on the strong solutions. Then the stability of weak solutions is verified by adapting the arguments for single-fluid equations \cite{AZ3,ZA1}. The key step in the asymptotic analysis is to show uniform-in-time bounds on the densities from above and below. Due to the complicated form of the pressure, classical methods used in \cite{HBV1,KN1,KZ,MY1,SV} cannot be applied here. However, thanks to the structure of the pressure, we are able to adapt the argument from \cite{Z1} so as to obtain the two-sided bounds; see Lemma \ref{lem9}. Based on these bounds, we  show the exponential decay of weak solution by choosing suitable test functions in the momentum equation and making another use of the structure of the pressure.

The functional spaces we use are standard. For brevity, we denote by $L^p$ the Lebesgue space $L^p((0,1))$ with the norm $\|\cdot\|_{L^p}$, and by  $H^1$ the Sobolev space $H^1((0,1))$.

Before stating our main results, we specify the meaning of weak solutions.
\begin{Definition}
Let $R_0,Q_0,u_0$ satisfy
\beq\label{bn16}
0<\underline{R_0}\leq R_0 \leq \overline{R_0}<\infty,
\eeq
\beq\label{bn17}
0<\underline{Q_0}\leq Q_0 \leq \overline{Q_0}<\infty,
\eeq
\beq\label{bn18}
u_0 \in L^2.
\eeq
A  triple  $(R,Q,u)$ is said to be a weak solution of (\ref{bn11})-(\ref{bn15}) on $[0,1]\times [0,T]$ provided that
\begin{itemize}
\item { $(\p_t R,\p_t Q,\p_y u)\in L^2(0,T;L^2)$,  \,\, $u\in L^{\infty}(0,T;L^2)$,}

\item {  $ 0< R(y,t),\,Q(y,t)<\infty$,  \, a.e. in  $(0,1)\times (0,T)$, }

\item { $\p_t \tau=\p_y u,\,R=R_0\tau_0 \tau^{-1},\,Q=Q_0\tau_0\tau^{-1}$, \, a.e. in $(0,1)\times (0,T)$,               }

\item { $\int_0^T\int_0^1 \left( u \p_t \phi-\left( \mu \f{\p_y u}{\tau}-Z^{\gamma_{+}}        \right)\p_y \phi \right)\dy\ds+\int_0^1 u_0\phi(\cdot,0)\dy=0,\\
     \text{for any } \phi\in C_c^{\infty}((0,1)\times [0,T))$.}
\end{itemize}
\end{Definition}
\begin{Remark}\label{rk}
Given $(R_0,Q_0)$ in (\ref{bn16})-(\ref{bn17}), we tacitly assume that $Z_0$ satisfies
\begin{equation}\label{bn22}
\left\{\begin{aligned}
& Q_0=\left(1-\f{R_0}{Z_0}\right) Z_0^{\gamma},\\
& R_0 \leq Z_0.\\
\end{aligned}\right.
\end{equation}
Clearly, the positive lower bound of $Z_0$ follows from (\ref{bn22})$_2$ and moreover $R_0<Z_0$ a.e. in $(0,1)$ in accordance with (\ref{bn22})$_1$. To get the upper bound of $Z_0$ we again make use of (\ref{bn22})$_1$. Indeed, suppose on the contrary that $Z_0>\max\{2\overline{R_0},(2\overline{Q_0})^{1/{\gamma}}\}$, then we would have
\[
\overline{Q_0}>Q_0=\left(1-\f{R_0}{Z_0}\right)Z_0^{\gamma}\geq \f{1}{2}Z_0^{\gamma}>\overline{Q_0},
\]
which is a contradiction. Therefore $Z_0$ must be bounded from above, more precisely
\[
Z_0\leq \max\left\{2\overline{R_0},(2\overline{Q_0})^{1/{\gamma}}\right\}.
\]
\end{Remark}

The main results of this paper are the following two theorems. The first one is concerned with the stability of weak solutions.
\begin{Theorem}\label{lsz1}
Let $\gamma_{\pm}>1$ and let \eqref{bn16}-\eqref{bn18} be satisfied.
%\beq\label{bn16}
%0<\underline{R_0}\leq R_0 \leq \overline{R_0}<\infty,
%\eeq
%\beq\label{bn17}
%0<\underline{Q_0}\leq Q_0 \leq \overline{Q_0}<\infty,
%\eeq
%\beq\label{bn18}
%u_0 \in L^2.
%\eeq
Then there exists a unique global-in-time weak solution to (\ref{bn11})-(\ref{bn15}). Moreover, if $(R,Q,u)$ and $(\widetilde{R},\widetilde{Q},\widetilde{u})$ are two weak solutions on $[0,1]\times [0,T]$ corresponding to the initial data $(R_0,Q_0,u_0)$ and $(\widetilde{R_0},\widetilde{Q_0},\widetilde{u_0})$, respectively, then
\[
\left( \|R-\widetilde{R}\|_{L^{\infty}(0,T;L^{\infty})}+ \|Q-\widetilde{Q}\|_{L^{\infty}(0,T;L^{\infty})}
+\|u-\widetilde{u}\|_{L^{\infty}(0,T;L^2)}+\|\p_y(u-\widetilde{u})\|_{L^2(0,T;L^2)} \right)
\]
\beq\label{bn19}
\leq C \left(\|R_0-\widetilde{R_0}\|_{L^{\infty}} +\|Q_0-\widetilde{Q_0}\|_{L^{\infty}}
+\|u_0-\widetilde{u_0}\|_{L^2} \right),
\eeq
where $C$ is a generic positive constant depending on $T$.
\end{Theorem}

The second theorem gives the large time behavior of weak solutions. More precisely, we show the asymptotic decay of weak solutions to $(R_{\infty},Q_{\infty},u_{\infty})$ --
 the unique steady state for problem (\ref{bn11})-(\ref{bn15}) given implicitly by
\begin{equation}\label{bn21}
\left\{\begin{aligned}
& Q_{\infty}=Q_0\tau_0 \tau_{\infty}^{-1},\,\,\,R_{\infty}=R_0\tau_0 \tau_{\infty}^{-1},\\
& \tau_{\infty}:=(R_{\infty}+Q_{\infty})^{-1},\\
& u_{\infty}=0,\,\,\,Z_{\infty}^{\gamma_{+}  } =C_{\star},\\
& Q_{\infty}=\left(1-\f{R_{\infty}}{Z_{\infty}}\right) Z_{\infty}^{\gamma},\,\,R_{\infty}\leq Z_{\infty},\\
& \int_0^1 \tau_{\infty}\dy=\int_0^1 \tau_0 \dy.
\end{aligned}\right.
\end{equation}
Here, $C_{\star}$ is the positive constant uniquely determined by $R_0,Q_0,\gamma_{\pm}$ and the conservation of mass (\ref{bn21})$_5$.

\begin{Theorem}\label{lsz2}
Let $(R,Q,u)$ be the unique weak solution to (\ref{bn11})-(\ref{bn15}) provided by Theorem \ref{lsz1}.
Then, for any $t\geq 0$, it holds
\beq\label{bn20}
\|(R-R_{\infty},Q-Q_{\infty},u-u_{\infty})\|_{L^2}  \leq C_1 \exp (-C_2 t),
\eeq
where $C_1$ and $C_2$ are generic positive constants independent of time.
\end{Theorem}
\begin{Remark}
Given suitably regular initial data, it can be shown, adapting the arguments from \cite{Z1,LYS1}, that
\[
\|(R-R_{\infty},Q-Q_{\infty},u-u_{\infty})\|_{H^1}  \leq C_1 \exp (-C_2 t).
\]
\end{Remark}

The rest of this paper is structured as follows. In Section \ref{sion21} we show global existence and uniqueness of strong solutions to (\ref{bn11})-(\ref{bn15}). In  Section  \ref{sion22} we prove the existence of global weak solutions via approximation based on regular solutions corresponding to regularized initial data and the weak convergence method. In Section \ref{sion3} we verify the stability of weak solutions. In Section \ref{sion4}, we obtain the exponential decay of weak solution to the unique steady state in $L^2$-norm with large initial data.

\section{Global existence of weak solutions}\label{sion2}
\subsection{Global well-posedness to (\ref{bn11})-(\ref{bn15})}\label{sion21}
In this subsection, we prove global existence and uniqueness of strong solutions to (\ref{bn11})-(\ref{bn15}) with large data. This will be useful in construction of weak solutions.
\begin{Proposition} \label {prop1}
Let (\ref{bn16})-(\ref{bn17}) be satisfied. Assume that
\beq\label{rdata}
(R_0,Q_0)\in H^1, \,u_0\in H_0^1.
\eeq

Then there exists a unique global strong solution $(R,Q,u)$ to (\ref{bn11})-(\ref{bn15}). Furthermore, for any $0<T<\infty$, it holds that
\beq\label{bm1}
C^{-1}\leq R(y,t),Q(y,t),Z(y,t) \leq C , \,\,\text{ for any }(y,t)\in [0,1]\times [0,T],
\eeq
\beq\label{bm2}
\|u\|_{ L^{\infty}(0,T;L^2)}+\|\p_y u\|_{L^2(0,T;L^2)  }  \leq C,
\eeq
\beq\label{bm3}
\|\p_t R\|_{L^2(0,T;L^2)  } +\|\p_t Q\|_{L^2(0,T;L^2)  } \leq C.
\eeq
\end{Proposition}

Local-in-time existence and uniqueness of strong solutions to (\ref{bn11})-(\ref{bn15}) is proved by the classical method based on the linearization of the problem and Banach fixed point theorem. We refer to \cite{EK1,NH} for similar calculations. Therefore, it remains to derive sufficient global a priori estimates so as to extend the local solution globally.

We start by giving the conservation of mass and the elementary energy inequality. To simplify the expression, we define
\beq\label{bm4}
\alpha:=\f{R}{Z}.
\eeq
In what follows, various positive constants are expressed by the same letter $C$ depending on $T$.

\begin{Lemma}\label{lem1}
Let $(R,Q,u)$ be a smooth solution to (\ref{bn11})-(\ref{bn15}) on $[0,1]\times [0,T]$ with regular initial data \eqref{rdata}, then we have
\beq\label{RQ}
R(y,t)=R_0\tau_0 \tau^{-1}(y,t),\quad Q(y,t)=Q_0\tau_0\tau^{-1}(y,t),\,\, \text{ for any } (y,t)\in [0,1]\times [0,T],
\eeq
\beq\label{bm5}
\int_0^1\tau(y,t)\dy=\int_0^1\tau_0(y)\dy,\,\, \text{ for any } t\in [0,T],
\eeq
\[
\sup_{0\leq t\leq T}\int_0^1\left(\f{1}{2} u^2 + \f{ (R_0\tau_0)^{\gamma_{+}}   }{\gamma_{+}-1} (\alpha \tau)^{-\gamma_{+}+1} +
  \f{ (Q_0\tau_0)^{\gamma_{-}}   }{\gamma_{-}-1} [(1-\alpha) \tau]^{-\gamma_{-}+1} \right)\dy
\]
\beq\label{bm6}
+\mu \int_0^T\int_0^1\f{(\p_y u)^2}{\tau}\dy\ds \leq C.
\eeq
\end{Lemma}
{\bf Proof.} The two first lines (\ref{RQ}) and (\ref{bm5}) follow from (\ref{bn11}) and (\ref{bn12}) immediately after integration with respect to time. To show (\ref{bm6}), we adopt to the Lagrangian coordinates the technique of pressure decomposition in \cite{BMZ}. One deduces from (\ref{bn4}), \eqref{notRQ}, and (\ref{bm4}) that
\beq\label{bm7}
\left( \f{R}{\alpha}  \right)^{\gamma_{+}}=\left( \f{Q}{1-\alpha}  \right)^{\gamma_{-}}.
\eeq
Thus we can decompose the pressure as
\beq\label{bm8}
Z^{ \gamma_{+}  }=\alpha \left( \f{R}{\alpha}  \right)^{\gamma_{+}}
+(1-\alpha)\left( \f{Q}{1-\alpha}  \right)^{\gamma_{-}}.
\eeq
Multiplying (\ref{bn13}) by $u$ and integrating by parts yields
\eq{\label{bm9}
0&=\f{1}{2}\f{d}{dt}\int_0^1 u^2 \dy -\int_0^1 Z^{ \gamma_{+}  } \p_y u \dy
+\mu \int_0^1 \f{(\p_y u)^2}{\tau}\dy\\
&=\f{1}{2}\f{d}{dt}\int_0^1 u^2 \dy -\int_0^1 Z^{ \gamma_{+}  } \p_t \tau \dy
+\mu \int_0^1 \f{(\p_y u)^2}{\tau}\dy.
}
The second term on the right-hand side of (\ref{bm9}) can be computed, with the help of (\ref{bm8}), through
\eq{\label{bm10}
-\int_0^1 Z^{ \gamma_{+}  } \p_t \tau \dy
&=-\int_0^1 \alpha \left( \f{R}{\alpha}  \right)^{\gamma_{+}} \p_t \tau \dy-\int_0^1
  (1-\alpha)\left( \f{Q}{1-\alpha}  \right)^{\gamma_{-}}     \p_t \tau \dy\\
&=-\int_0^1 \left( \f{R_0 \tau_0}{\alpha \tau}  \right)^{\gamma_{+}} \alpha  \p_t \tau \dy
-\int_0^1
\left( \f{Q_0\tau_0}{(1-\alpha)\tau}  \right)^{\gamma_{-}}    (1-\alpha)   \p_t \tau \dy\\
&=-\int_0^1 \left( \f{R_0 \tau_0}{\alpha \tau}  \right)^{\gamma_{+}} \p_t (\alpha   \tau ) \dy
+\int_0^1 \left( \f{R_0 \tau_0}{\alpha \tau}  \right)^{\gamma_{+}} \tau  \p_t \alpha \dy\\
&\quad-\int_0^1
\left( \f{Q_0\tau_0}{(1-\alpha)\tau}  \right)^{\gamma_{-}}  \p_t  [(1-\alpha) \tau ]\dy
+\int_0^1
\left( \f{Q_0\tau_0}{(1-\alpha)\tau}  \right)^{\gamma_{-}}  \tau  \p_t (1-\alpha)  \dy\\
&=-\int_0^1 \left( \f{R_0 \tau_0}{\alpha \tau}  \right)^{\gamma_{+}} \p_t (\alpha   \tau ) \dy
-\int_0^1
\left( \f{Q_0\tau_0}{(1-\alpha)\tau}  \right)^{\gamma_{-}}  \p_t  [(1-\alpha) \tau ]\dy\\
&=\f{d}{dt} \int_0^1 \f{ (R_0\tau_0)^{\gamma_{+}}   }{\gamma_{+}-1} (\alpha \tau)^{-\gamma_{+}+1}
\dy +\f{d}{dt} \int_0^1 \f{ (Q_0\tau_0)^{\gamma_{-}}   }{\gamma_{-}-1} [(1-\alpha) \tau]^{-\gamma_{-}+1} \dy,}
where we have used \eqref{RQ} for the second equality, and (\ref{bm7}) in the fourth equality. Thus, combining (\ref{bm9}) and (\ref{bm10}) gives rise to (\ref{bm6}). This finishes the proof of Lemma \ref{lem1}. $\Box$

By virtue of Lemma \ref{lem1} and the specific mathematical structure of the equations, we are able to show the upper and lower bounds for $R$ and $Q$. This plays a crucial role in the proof of Proposition \ref{prop1}. The idea of proof comes from \cite{AKS}.

\begin{Lemma}\label{lem2}
Let the assumptions of Lemma \ref{lem1} be satisfied. Then, there exists a constant $C>0$ such that
\beq\label{bm11}
C^{-1} \leq R(y,t),Q(y,t),Z(y,t)\leq C, \,\, \text{ for any }(y,t)\in [0,1]\times [0,T].
\eeq
\end{Lemma}
{\bf Proof.} We observe first that the positiveness of $R$ and $Q$ follow from the method of characteristics and the regular initial data. Given positive $R$ and $Q$, there exists a unique $Z$ satisfying (\ref{bn9})-(\ref{bn10}). This fact can be justified easily and we refer to Lemma 2.1 in \cite{BMZ} for the details. That is, $Z$ can be regarded as a function of $R$ and $Q$. Moreover, due to \eqref{RQ}  and \eqref{bn9}, $Z$ can be also seen  as a function of $y$ and $\tau$.

In accordance with (\ref{bm5}) and due to continuity of $\tau$, for any $t\in [0,T]$, there exists $a(t)\in [0,1]$ such that
\[
\tau(a(t),t)=\int_0^1\tau_0(y)\dy.
\]
Integrating (\ref{bn13}) with respect to time over $(0,t)$, and  then with respect to space over $(a(t),y)$, followed by taking exponentials on both sides of the resulting equation, we arrive at
\beq\label{bm12}
Y(t)\tau(y,t)=\tau_0(y)B(y,t)\exp \left( \f{1}{\mu} \int_0^t Z^{\gamma_{+}}(y,\tau)(y,s)\ds \right),
\eeq
where
\[
Y(t):= \tau_0(a(t)) \tau^{-1}(a(t),t) \exp \left(   \f{1}{\mu} \int_0^t Z^{\gamma_{+}}(a(t),s)\ds       \right),
\]
\[
B(y,t):= \exp \left(  \f{1}{\mu}\int_{a(t)}^y [u(\xi,t)-u_0(\xi)] d\xi     \right).
\]
Next, in order to get the two-sided bounds from the representation formula (\ref{bm12}), one needs the bounds for $Y$ and $B$ from above and below. Obviously, by (\ref{bm6}) and Cauchy-Schwarz's inequality, it follows that
\beq\label{bm13}
C^{-1} \leq B(y,t)\leq C, \text{ for any }(y,t)\in [0,1]\times [0,T].
\eeq
It is also clear that
\beq\label{bm14}
C^{-1} \leq Y(t), \text{ for any }t\in [0,T].
\eeq
For the purpose of obtaining the upper bound of $Y$, we rewrite (\ref{bm12}) in another form. Noticing that
\[
\f{\p}{\p t} \exp \left( \f{1}{\mu} \int_0^t Z^{\gamma_{+}}(y,\tau)(y,s)\ds \right)
\]
\[
=\f{1}{\mu}Z^{\gamma_{+}}(y,\tau) \exp \left( \f{1}{\mu} \int_0^t Z^{\gamma_{+}}(y,\tau)(y,s)\ds \right)
\]
\[
=\f{1}{\mu}Z^{\gamma_{+}}(y,\tau) Y(t)\tau (y,t) \tau_0^{-1}(y)B^{-1}(y,t),
\]
and integrating it with respect to time, we get
\[
\exp \left( \f{1}{\mu} \int_0^t Z^{\gamma_{+}}(y,\tau)(y,s)\ds \right)
\]
\[
=1+\f{1}{\mu} \int_0^t Z^{\gamma_{+}}(y,\tau)(y,s) \tau(y,s) Y(s) \tau_0^{-1}(y) B^{-1}(y,s) \ds.
\]
Inserting the above relation back to (\ref{bm12}) yields
\[
Y(t)\tau(y,t)
\]
\beq\label{bm15}
=\tau_0(y)B(y,t) \left(   1+\f{1}{\mu} \int_0^t Z^{\gamma_{+}}(y,\tau)(y,s) \tau(y,s) Y(s) \tau_0^{-1}(y) B^{-1}(y,s) \ds   \right).
\eeq
In light of (\ref{bm8}), we have
\[
\int_0^1 \tau Z^{\gamma_{+}}  \dy
=\int_0^1 \left[ \tau \alpha \left( \f{R}{\alpha}  \right)^{\gamma_{+}}
+\tau(1-\alpha)\left( \f{Q}{1-\alpha}  \right)^{\gamma_{-}} \right] \dy
\]
\[
= \int_0^1 \left(  (R_0\tau_0)^{\gamma_{+}}    (\alpha \tau)^{-\gamma_{+}+1} +
 (Q_0\tau_0)^{\gamma_{-}}    [(1-\alpha) \tau]^{-\gamma_{-}+1}  \right) \dy
\]
\beq\label{bm16}
\leq C,
\eeq
due to (\ref{bm6}). Now we integrate both sides of (\ref{bm15})   with respect to space over $(0,1)$, and make use of (\ref{bm5}), (\ref{bm13}) and (\ref{bm16}) to conclude that
\[
Y(t)\leq C+ C \int_0^t Y(s)\ds.
\]
A straightforward application of Gronwall's inequality gives
\beq\label{bm17}
Y(t)\leq C,\,\, \text{ for any }t\in [0,T].
\eeq
And so, using (\ref{bm13}) and (\ref{bm17}), we deduce from (\ref{bm12}) that
\beq\label{bm18}
C^{-1}\leq \tau(y,t), \,\,\text{ for any }(y,t)\in [0,1]\times [0,T].
\eeq
It remains to get the upper bound of $\tau$. To this end, we rewrite (\ref{bn9}) as
\beq\label{bm19}
Q_0\tau_0\tau^{-1}=\left( 1-\f{R_0\tau_0\tau^{-1}}{Z} \right)Z^{\gamma}.
\eeq
Differentiating both sides of (\ref{bm19})  with respect to $\tau$ leads to
\[
-Q_0\tau_0\tau^{-2} =\gamma Z^{\gamma-1} \f{\p Z}{\p \tau}
-\left( R_0\tau_0\tau^{-1}(\gamma-1) Z^{\gamma-2} \f{\p Z}{\p \tau}- R_0\tau_0\tau^{-2}Z^{\gamma-1}      \right),
\]
or equivalently,
\beq\label{bm20}
\f{\p Z}{\p \tau}=-\f{Q_0\tau_0\tau^{-2}+ R_0\tau_0\tau^{-2}Z^{\gamma-1}    }{\gamma Z^{\gamma-1}-R_0\tau_0\tau^{-1}(\gamma-1) Z^{\gamma-2}}.
\eeq
The denominator is positive as we have
\eq{\label{bm21}
\gamma Z^{\gamma-1}-R_0\tau_0\tau^{-1}(\gamma-1) Z^{\gamma-2}
= Z^{\gamma-2} [ \gamma Z- (\gamma-1)R  ]\\
= Z^{\gamma-2} [ \gamma (Z-R) +R]\geq Z^{\gamma-2}R>0,
}
due to (\ref{bn10}), assumption about the smoothness of the solution, and Remark \ref{rk}. Therefore,
\beq\label{bm22}
\f{\p Z^{\gamma_{+}}}{\p \tau}=\gamma_{+}Z^{\gamma_{+}-1}\f{\p Z}{\p \tau}<0
\eeq
which means that the pressure is decreasing with respect to $\tau$. From (\ref{bm12}) and the two-sided bounds of $Y(t), B(y,t)$, we see
 \[
 \tau \leq C \exp \left(  \int_0^t Z^{\gamma_{+}}(y,\tau)(y,s)\ds \right)
 \]
 \[
 \leq C \exp \left(  \int_0^t Z^{\gamma_{+}}(y,\underline{\tau})(y,s)\ds \right),
 \]
 since $\tau \geq \underline{\tau}$ due to (\ref{bm18}) and the fact that the pressure is decreasing with respect to $\tau$. Notice that $Z(y,\underline{\tau})$ obeys
 \[
 Q_0\tau_0 (\underline{\tau})^{-1}=\left( 1-\f{R_0\tau_0 (\underline{\tau})^{-1}}{Z(y,\underline{\tau})} \right)Z(y,\underline{\tau})^{\gamma}.
 \]
 Arguing as in Remark \ref{rk} one sees that $Z(y,\underline{\tau})$ must be bounded from above with upper bound depending only on $R_0,Q_0,\underline{\tau}$. Thus,
\beq\label{bm23}
 \tau(y,t)\leq C, \,\,\text{ for any }(y,t)\in [0,1]\times [0,T].
\eeq
As an immediate consequence of (\ref{bm18}), (\ref{bm23}) and the relations \eqref{RQ}, we finally verify (\ref{bm11}) for $R$ and $Q$. Notice that the two-sided bounds for $Z$ follow from the same argument as in Remark \ref{rk}. This completes the proof of Lemma \ref{lem2}. $\Box$

We conclude directly from Lemmas \ref{lem1}-\ref{lem2} that
\begin{Corollary}\label{lem3} Under the assumptions of Lemma \ref{lem1}, we have
\beq\label{bm24}
\|u\|_{ L^{\infty}(0,T;L^2)}+\|\p_y u\|_{L^2(0,T;L^2)  }  \leq C,
\eeq
\beq\label{bm25}
\|\p_t R\|_{L^2(0,T;L^2)  } +\|\p_t Q\|_{L^2(0,T;L^2)  } \leq C.
\eeq
\end{Corollary}

We are now in a position to prove that $(R,Q,u)$ are more regular in order to finish the proof of Proposition \ref{prop1}. This process is quite standard due to Lemma \ref{lem2}. Thus the estimates are listed below and the details are omitted here and we refer to \cite{AKS} for similar calculations.
\begin{Lemma}\label{lem4}
The strong solution to (\ref{bn11})-(\ref{bn15}) on $[0,1]\times [0,T]$ satisfies
\beq\label{bm26}
\|\p_y R\|_{ L^{\infty}(0,T;L^2)}+\|\p_y Q\|_{ L^{\infty}(0,T;L^2)}  \leq C,
\eeq
\beq\label{bm27}
\|\p_y u\|_{ L^{\infty}(0,T;L^2)} +\|\p_{yy} u\|_{L^2(0,T;L^2)  } \leq C.
\eeq
\end{Lemma}

With Lemmas \ref{lem1}-\ref{lem4} and the Corollary \ref{lem3} at hand, the local-in-time solution can be extended globally in a routine manner. Uniqueness of solutions is proved by the classical energy method. This finishes the proof of Proposition \ref{prop1}. $\Box$

\subsection{Existence of weak solutions}\label{sion22}
The main task of this subsection is to construct global-in-time weak solutions to (\ref{bn11})-(\ref{bn15}) using approximation based on regular solutions. We start from regularizing the initial data $(R_0,Q_0,u_0)$ in such a way that $\{(R_0^{\ep},Q_0^{\ep},u_0^{\ep})\}_{\ep>0}$ satisfy
\[
(R_0^{\ep},Q_0^{\ep}) \in C^2 ([0,1]),\,\,C^{-1} \leq R_0^{\ep},Q_0^{\ep} \leq C,\,\, u_0^{\ep}\in C_c^2 ((0,1)),
\]
\[
(R_0^{\ep},Q_0^{\ep},u_0^{\ep}) \rightarrow (R_0,Q_0,u_0) \text{ strongly in }L^2 \text{ as }\ep \rightarrow 0.
\]
Moreover, we define
\[
\tau_0^{\ep}:=(R_0^{\ep}+Q_0^{\ep})^{-1},\,\,\tau_{\ep}:=(R_{\ep}+Q_{\ep})^{-1}.
\]
Therefore, it follows from Proposition \ref{prop1} that there exists a unique global strong solution $(R_{\ep},Q_{\ep},u_{\ep})$ to (\ref{bn11})-(\ref{bn15}) with initial data $(R_0^{\ep},Q_0^{\ep},u_0^{\ep})$. Furthermore, from Proposition \ref{prop1} we conclude the following uniform-in-$\ep$ estimates:
\beq\label{bm28}
C^{-1}\leq R_{\ep}(y,t),Q_{\ep}(y,t) \leq C , \,\,\text{ for any }(y,t)\in [0,1]\times [0,T],
\eeq
\beq\label{bm29}
\|u_{\ep}\|_{ L^{\infty}(0,T;L^2)}+\|\p_y u_{\ep}\|_{L^2(0,T;L^2)  }  \leq C,
\eeq
\beq\label{bm30}
\|\p_t R_{\ep}\|_{L^2(0,T;L^2)  } +\|\p_t Q_{\ep}\|_{L^2(0,T;L^2)  } \leq C.
\eeq
From (\ref{bm28})-(\ref{bm30}) it follows that there exists a subsequence of $\{(R_{\ep},Q_{\ep},u_{\ep})\}_{\ep>0}$, not relabeling, such that as $\ep\rightarrow 0$,
\beq\label{bm31}
(R_{\ep},Q_{\ep})\rightarrow (R,Q) \text{ weakly}-\ast \text{ in } L^{\infty}(0,T;L^{\infty}),
\eeq
\beq\label{bm32}
u_{\ep}\rightarrow u  \text{ weakly}-\ast \text{ in }L^{\infty}(0,T;L^2),
\eeq
\beq\label{bm33}
(\p_t R_{\ep},\p_t Q_{\ep},\p_y u_{\ep})\rightarrow (\p_t R,\p_t Q,\p_y u)\text{ weakly in }L^2(0,T;L^2),
\eeq
for some limit triple $(R,Q,u)$, and we infer from (\ref{bm28})-(\ref{bm33}) that
\beq\label{bm34}
C^{-1}\leq R(y,t),Q(y,t) \leq C , \,\,\text{ for a.e. }(y,t)\in (0,1)\times (0,T),
\eeq
\beq\label{bm35}
\|u\|_{ L^{\infty}(0,T;L^2)}+\|\p_y u\|_{L^2(0,T;L^2)  }  \leq C,
\eeq
\beq\label{bm36}
\|\p_t R\|_{L^2(0,T;L^2)  } +\|\p_t Q\|_{L^2(0,T;L^2)  } \leq C.
\eeq

The weak convergence results (\ref{bm31})-(\ref{bm33}) are not sufficient to pass to the limit in (\ref{bn11})-(\ref{bn15}), in particular, in the strongly nonlinear  pressure function. For the moment we only know that $Z_{\ep}$ is the unique solution of
\[
Q_{\ep}=\left( 1-\f{R_{\ep}}{Z_{\ep}} \right)Z_{\ep}^{\gamma} ,\,\, R_{\ep}\leq Z_{\ep}.
\]
To identify the pressure term, it suffices to verify that the pointwise limit of $\{Z_{\ep}\}_{\ep>0}$ is the unique solution of
\[
Q=\left(1-\f{R}{Z}\right) Z^{\gamma},\,\,R \leq Z,
\]
for which we need the strong convergence of the sequence $\{Z_\ep\}_{\ep>0}$.
In fact, since $Q_{\ep}=Q_0^{\ep}\tau_0^{\ep}\tau_{\ep}^{-1}$, $R_{\ep}=R_0^{\ep}\tau_0^{\ep}\tau_{\ep}^{-1}$,  $Z_{\ep}$ can also be regarded as $Z_{\ep}=Z_{\ep}(Q_0^{\ep},R_0^{\ep},\tau_0^{\ep},\tau_{\ep})$. Therefore the strong convergence of $\{Z_\ep\}_{\ep>0}$ will follow from that of $\{\tau_{\ep}\}_{\ep>0}$. The necessary compactness property in space is provided by the following lemma.

\begin{Lemma}\label{lem5}
For any $0<h<1$, there holds
\beq\label{bm37}
\|\Delta_h \tau_{\varepsilon}\|_{L^{\infty}(0,T;L^2)  }
\leq C ( \|\Delta_h R_0\|_{L^2}+\|\Delta_h Q_0\|_{L^2} +h ),
\eeq
where $\Delta_h F(y):=F(y+h)-F(y)$ is the translation in spatial variable with the step $h$.
\end{Lemma}
{\bf Proof.} Similarly to Lemma \ref{lem2}, we need a representation formula of $\tau_{\ep}$. By setting
\[
\sigma_{\ep}:= \mu \f{\p_y u_{\ep}}{\tau_{\ep}}-Z_{\ep}^{\gamma_{+}},
\]
and recalling that $(R_{\ep},Q_{\ep},u_{\ep})$ solves (\ref{bn11})-(\ref{bn13}) in the strong sense, it follows that
\[
\p_t \tau_{\ep}= \f{1}{\mu} \tau_{\ep}(\sigma_{\ep}+Z_{\ep}^{\gamma_{+}}).
\]
Multiplying the above identity both sides by $\exp\left(-\f{1}{\mu}\int_0^t\sigma_{\ep}(y,s)\ds\right)$ and integrating over time, one gets the relation
\beq\label{bm38}
\tau_{\ep}=D_{\ep}\left(\tau^{\ep}_0+\f{1}{\mu}
\int_0^t D_{\ep}^{-1}(y,s)
\left(\tau_{\ep}Z_{\ep}^{\gamma_{+}}
 \right)(y,s)\ds\right),
\eeq
where
\[
D_{\ep}(y,t):=\exp\left(\f{1}{\mu}\int_0^t\sigma_{\ep}(y,s)\ds\right).
\]
By definition it holds that
\[
\Delta_h\tau_{\ep}(y,t)= \Delta_h D_{\ep}(y,t)\left(\tau^{\ep}_0(y+h)+\f{1}{\mu}
\int_0^t D_{\ep}^{-1}
\left(\tau_{\ep}Z_{\ep}^{\gamma_{+}}\right)(y+h,s)\ds\right)
\]
\[
+  D_{\ep}(y,t)\left(\Delta_h\tau^{\ep}_0(y)+\f{1}{\mu}\int_0^t D_{\ep}^{-1}(y+h,s)\Delta_h\left(\tau_{\ep}Z_{\ep}^{\gamma_{+}}\right)
(y,s) \ds\right)
\]
\beq\label{bm39}
-  D_{\ep}(y,t)\left( \f{1}{\mu} \int_0^t D_{\ep}^{-1}(y+h,s)D_{\ep}^{-1}(y,s)
\left(\tau_{\ep}Z_{\ep}^{\gamma_{+}}\right)(y,s)\Delta_h D_{\ep}(y,s) \ds\right).
\eeq
Thanks to (\ref{bm28}), we have
\beq\label{bm40}
C^{-1}\leq Z_{\ep}(y,t),\,\,D_{\ep}(y,t) \leq C , \,\,\text{ for any }(y,t)\in [0,1]\times [0,T].
\eeq
Indeed, the lower bound of $Z_{\ep}$ follows readily from (\ref{bm28}) and the relation $R_{\ep}\leq Z_{\ep}$; the upper bound of $Z_{\ep}$ is verified by the relation
\[
Q_{\ep}=\left( 1-\f{R_{\ep}}{Z_{\ep}} \right)Z_{\ep}^{\gamma},
\]
and the two-sided bounds of $R_{\ep},Q_{\ep}$ through the argument as in Remark \ref{rk}. The delicate issue is to compute $\Delta_h\left(\tau_{\ep}Z_{\ep}^{\gamma_{+}}\right)$. In fact,
\beq\label{bm41}
\Delta_h\left(\tau_{\ep}Z_{\ep}^{\gamma_{+}}\right)=Z_{\ep}^{\gamma_{+}} \Delta_h \tau_{\ep} +
\tau_{\ep}(y+h,t)\Delta_h Z_{\ep}^{\gamma_{+}}.
\eeq
Furthermore, as pointed out before, $Z_{\ep}$ can be regarded as a function $Z_{\ep}=Z_{\ep}(Q_0^{\ep},R_0^{\ep},\tau_0^{\ep},\tau_{\ep})$. Thus by the mean value theorem
\[
\Delta_h Z_{\ep}=\f{\p Z_{\ep}}{\p Q_0^{\ep}} \Delta_h  Q_0^{\ep}+
\f{\p Z_{\ep}}{\p R_0^{\ep}} \Delta_h  R_0^{\ep}+
\f{\p Z_{\ep}}{\p \tau_0^{\ep}} \Delta_h  \tau_0^{\ep}+
\f{\p Z_{\ep}}{\p \tau_{\ep}} \Delta_h  \tau_{\ep}.
\]
Subsequent differentiations of \eqref{bm19} with respect to $Q_0^\ep, R_0^\ep, \tau_0^\ep$, and $\tau_\ep$ give rise to
\[
\f{\p Z_{\ep}}{\p Q_0^{\ep}}=\f{\tau_0^{\ep}\tau_{\ep}^{-1}}{\gamma Z_{\ep}^{\gamma-1}-R_0^{\ep}\tau_0^{\ep}\tau_{\ep}^{-1}(\gamma-1) Z_{\ep}^{\gamma-2}},
\]
\[
\f{\p Z_{\ep}}{\p R_0^{\ep}}=\f{\tau_0^{\ep}\tau_{\ep}^{-1}Z_{\ep}^{\gamma-1}}{\gamma Z_{\ep}^{\gamma-1}-R_0^{\ep}\tau_0^{\ep}\tau_{\ep}^{-1}(\gamma-1) Z_{\ep}^{\gamma-2}},
\]
\[
\f{\p Z_{\ep}}{\p\tau_0^{\ep}}=\f{Q_0^{\ep}\tau_{\ep}^{-1}+R_0^{\ep}\tau_{\ep}^{-1}Z_{\ep}^{\gamma-1}}
{\gamma Z_{\ep}^{\gamma-1}-R_0^{\ep}\tau_0^{\ep}\tau_{\ep}^{-1}(\gamma-1) Z_{\ep}^{\gamma-2}},
\]
\[
\f{\p Z_{\ep}}{\p \tau_{\ep}}=-\f{Q_0^{\ep}\tau_0^{\ep}\tau_{\ep}^{-2}+ R_0^{\ep}\tau_0^{\ep}\tau_{\ep}^{-2}Z_{\ep}^{\gamma-1}    }{\gamma Z_{\ep}^{\gamma-1}-R_0^{\ep}\tau_0^{\ep}\tau_{\ep}^{-1}(\gamma-1) Z_{\ep}^{\gamma-2}}.
\]
In view of (\ref{bm21}), (\ref{bm28}) and (\ref{bm40}), it follows that
\beq\label{bm42}
\left\|\left(\f{\p Z_{\ep}}{\p Q_0^{\ep}},\f{\p Z_{\ep}}{\p R_0^{\ep}},\f{\p Z_{\ep}}{\p\tau_0^{\ep}},\f{\p Z_{\ep}}{\p \tau_{\ep}}\right)\right\|_{L^{\infty}(0,T;L^{\infty})}\leq C.
\eeq
As a consequence, we conclude from (\ref{bm28})-(\ref{bm29}) and (\ref{bm39})-(\ref{bm42}) that
\[
\|\Delta_h\tau_{\ep}(y,t)\|_{L^2} \leq  C\Big(\|\Delta_h D_{\ep}(y,t)\|_{L^2}+
 \|\Delta_h\tau^{\ep}_0\|_{L^2}
\]
\[
+\int_0^t \left( \|\Delta_h D_{\ep}(y,s)\|_{L^2}
 +\|\Delta_h\tau_{\ep}(y,s)\|_{L^2}
 +\|\Delta_h Q^{\ep}_0\|_{L^2}+\|\Delta_h R^{\ep}_0\|_{L^2} \right)\ds \Big)
\]
\[
 \leq  C\left(h\|u_{\ep}-u_0^{\ep}\|_{L^\infty(0,T;L^2)}+
 \|\Delta_h Q^{\ep}_0\|_{L^2}
  + \|\Delta_h R^{\ep}_0\|_{L^2}
  +\int_0^t\|\Delta_h\tau_{\ep}(y,s)\|_{L^2} \ds \right)
\]
\beq\label{bm43}
 \leq C \left(h+\|\Delta_h Q_0\|_{L^2}+\|\Delta_h R_0\|_{L^2}+ \int_0^t\|\Delta_h\tau_{\ep}(y,s)\|_{L^2} \ds \right).
\eeq
Finally, (\ref{bm37}) follows from (\ref{bm43}) immediately by invoking Gronwall's inequality. The proof of Lemma \ref{lem5} is thus finished. $\Box$

Based on Lemma \ref{lem5}, the relation $\p_t \tau_{\ep}=\p_y u_{\ep}$ and (\ref{bm29}), we see
\[
\|\tau_{\ep}(\cdot+h,\cdot+s)-\tau_{\ep}\|_{L^{\infty}(0,T-s;L^2)}
\leq C\left(\|\Delta_h Q_0\|_{L^2}+\|\Delta_h R_0\|_{L^2}+h+s^{\f{1}{2}}\right),
\]
for any $0<h<1$, $0<s<T$. This particularly implies the strong convergence of $\{\tau_{\ep}\}_{\ep>0}$ to $\tau$ in $L^2(0,T;L^2)$ and furthermore in $L^p(0,T;L^p)$ for any $1\leq p<\infty$. Consequently, it holds that
\[
Q_{\ep}\rightarrow Q_0\tau_0\tau^{-1},\,\, R_{\ep}\rightarrow R_0\tau_0\tau^{-1},\,\, \text{ a.e. in }(0,1)\times (0,T),
\]
which yields
\beq\label{bm44}
Q=Q_0\tau_0 \tau^{-1},\,\,R=R_0\tau_0 \tau^{-1}.
\eeq
Recalling that $Z_{\ep}=Z_{\ep}(Q_0^{\ep},R_0^{\ep},\tau_0^{\ep},\tau_{\ep})$, we find $Z_{\ep}$ converges to some limit function $Z$ almost everywhere. Upon passing to the limit in the relations
\[
Q_{\ep}=\left( 1-\f{R_{\ep}}{Z_{\ep}} \right)Z_{\ep}^{\gamma} ,\,\, R_{\ep}\leq Z_{\ep},
\]
we conclude from (\ref{bm40}) that $\{Z_{\ep}\}_{\ep>0}$ converges to $Z$ strongly in $L^p(0,T;L^p)$ for any $1\leq p<\infty$ and $Z$ solves exactly
\beq\label{bm45}
Q=\left(1-\f{R}{Z}\right) Z^{\gamma},\,\,R \leq Z.
\eeq
This finishes the proof of existence of a weak solution.
\section{Stability of weak solutions}\label{sion3}
In the present section, we show Lipschitz continuous dependence on the initial data of weak solutions, i.e., we prove our first main Theorem \ref{lsz1}. We remark that the proof relies on the structure of the equations. As a preliminary step, we state the following lemma, the proof of which is omitted as it is similar to relevant results from  \cite{LYS1,ZA1}.
\begin{Lemma}\label{lem6}
Let $(R,Q,u)$ be a weak solution to (\ref{bn11})-(\ref{bn15}). Then
\[
\tau(y,t) = \exp\left(\f{1}{\mu}\int_0^t\sigma(y,s)\ds\right)
\]
\beq\label{bv1}
          \times  \left( \tau_0+\f{1}{\mu}\int_0^t\exp\left(-\f{1}{\mu}\int_0^{\xi}\sigma(y,s)\ds\right)
          \left(\tau Z^{\gamma_{+}}\right)(y,\xi) d\xi \right),
\eeq
and
\beq\label{bv2}
\int_0^t\sigma(y,s)\ds=(\mathcal{I}(u-u_0))(y,t)+\int_0^t\langle\sigma(\cdot,s)\rangle \ds,
\eeq
where
\[
\sigma(y,t):=\left( \mu \f{\p_y u}{\tau}-Z^{\gamma_{+}}\right)(y,t),
\]
\[
\mathcal{I}f(y):=\int_0^y f(\xi)d\xi-\left<\int_0^y f(\xi)d\xi\right>,\,\, \langle f\rangle:=\int_0^1 f(y)\dy.
\]
\end{Lemma}

  To verify the stability estimate (\ref{bn19}) from Theorem \ref{lsz1}, we follow the arguments in \cite{AZ3,LYS1}. Let us start from introducing the following notation:
\begin{equation*}
\left\{\begin{aligned}
& (\Delta \tau,\Delta R,\Delta Q, \Delta u):=(\tau-\widetilde{\tau},R-\widetilde{R},Q-\widetilde{Q},u-\widetilde{u}), \\
& (\Delta \tau_0,\Delta R_0,\Delta Q_0, \Delta u_0):=(\tau_0-\widetilde{\tau_0},R_0-\widetilde{R_0},Q_0-\widetilde{Q_0},u_0-\widetilde{u_0}),\\
& \Delta\sigma:=\sigma-\widetilde{\sigma},\,\,
\widetilde{\sigma}:=\mu\f{\p_y \widetilde{u}}{\widetilde{\tau}}
-(\widetilde{Z})^{\gamma_{+}}, \\
& D:=\exp\left(\f{1}{\mu}\int_0^t\sigma(y,s)\ds\right),\,\,
\widetilde{D}:=\exp\left(\f{1}{\mu}\int_0^t\widetilde{\sigma}(y,s)\ds\right),\\
& \widetilde{\vr}:=\widetilde{R}+\widetilde{Q},\,\,\Delta \vr:=\vr-\widetilde{\vr}.
\end{aligned}\right.
\end{equation*}
Recalling that $Q=Q_0\tau_0 \tau^{-1},\,R=R_0\tau_0 \tau^{-1}$, one has in light of  uniform bounds  for $R,Q$ from below and from above, i.e., (\ref{bm34}), that
\eq{\label{bv3}
&|\Delta R| \leq C  \left(|\Delta R_0| +|\Delta \tau_0|+|\Delta \tau|            \right)
 \leq C  \left(|\Delta R_0| +|\Delta Q_0|+|\Delta \tau|            \right);\\
&|\Delta Q| \leq C  \left(|\Delta Q_0| +|\Delta \tau_0|+|\Delta \tau|            \right) \leq C  \left(|\Delta R_0| +|\Delta Q_0|+|\Delta \tau|            \right).
}
Consequently, in order to estimate $L^{\infty}(0,T;L^{\infty})$-norm of $\Delta R$ and $\Delta Q$, it suffices to control $\|\Delta \tau\|_{L^{\infty}(0,T;L^{\infty})}$. This is the key step in proving stability of weak solutions. We follow the idea in \cite{AZ3,ZA1} to accomplish this goal.
\begin{Lemma}\label{lem7}
Let the assumptions of Theorem \ref{lsz1} be fulfilled, then we have
\eq{\label{bv4}
&\|\Delta\tau\|_{L^{\infty}(0,t;L^{\infty})} \leq C \Big(\|\Delta R_0\|_{L^{\infty}}+\|\Delta Q_0\|
_{L^{\infty}}+\|\Delta u_0\|_{L^2}\\
&\hspace{4cm}+\|\Delta u\|_{L^{\infty}(0,t;L^2)}
+\|\p_y (\Delta u)\|_{L^2(0,t;L^2)}\Big)
}
for any  $t\in (0, T]$, where $C$ denotes generic positive constant depending on $T$.
\end{Lemma}
{\bf Proof.} It follows from (\ref{bv1}) that
\[
\Delta\tau  = D\left\{\Delta \tau_0+\f{1}{\mu}\int_0^t\left( \tau Z^{\gamma_{+}} \left(D^{-1}-(\widetilde{D})^{-1}\right)+\f{\tau Z^{\gamma_{+}}-\widetilde{\tau} (\widetilde{Z})^{\gamma_{+}}  }
{\widetilde{D}}\right)\ds\right\}
\]
\beq\label{bv5}
+ \left(D-\widetilde{D}\right)\left(\widetilde{\tau_0}+\f{1}{\mu}\int_0^t\f{\widetilde{\tau} (\widetilde{Z})^{\gamma_{+}}}{\widetilde{D}}\ds\right).
\eeq
Similarly to (\ref{bm40}), it holds that
\beq\label{bv6}
C^{-1}\leq \left(Z,\widetilde{Z},D,\widetilde{D}\right)(y,t) \leq C , \,\,\text{ for a.e. }(y,t)\in (0,1)\times (0,T).
\eeq
Indeed, the upper bound of $Z$ and $\widetilde{Z}$ is derived by the same argument as (\ref{bm40}). Based on (\ref{bm34}) and (\ref{bv6}), we observe that
\eq{\label{bv7}
\left|\tau Z^{\gamma_{+}}-\widetilde{\tau} (\widetilde{Z})^{\gamma_{+}} \right|
\leq C &\Big( |\Delta \tau|+ \left\|\f{\p Z}{\p \tau}\right\|_{L^{\infty}(0,T;L^{\infty})}|\Delta \tau| +\left\|\f{\p Z}{\p Q_0}\right\|_{L^{\infty}(0,T;L^{\infty})} |\Delta Q_0|\\
&\hspace{1cm}+ \left\|\f{\p Z}{\p R_0}\right\|_{L^{\infty}(0,T;L^{\infty})} |\Delta R_0| +\left\|\f{\p Z}{\p \tau_0}\right\|_{L^{\infty}(0,T;L^{\infty})}|\Delta \tau_0| \Big)\\
\leq C& \left( |\Delta \tau|+|\Delta R_0|+|\Delta Q_0|  \right),
}
where we used a version of \eqref{bm42} for the limit functions. Therefore, we deduce from (\ref{bv5})-(\ref{bv7}) that
\[
|\Delta\tau|\leq C\left(|\Delta R_0|+|\Delta Q_0|+\int_0^t\left(\left|\int_0^{s}\Delta\sigma d\xi\right|+|\Delta\tau|\right) \ds\right)
+C\left|\int_0^t\Delta\sigma\ds\right|;
\]
whence
\beq\label{bv8}
\|\Delta\tau(\cdot,t)\|_{L^{\infty}}\leq C\left(\|\Delta R_0\|_{L^{\infty}}
+\|\Delta Q_0\|_{L^{\infty}}+\left\|\int_0^{s}\Delta\sigma d\xi\right\|_{L^{\infty}
(0,t;L^{\infty})}+\int_0^t\|\Delta\tau(\cdot,s)\|_{L^{\infty}} \ds \right).
\eeq
The rest of the proof follows the same lines as \cite{LYS1}, and we write down the details only for the convenience of the reader. First, from the identity (\ref{bv2}) and H\"{o}lder's inequality we obtain
\beq\label{bv9}
\left\|\int_0^{s}\Delta\sigma d\xi\right\|_{L^{\infty}
(0,t;L^{\infty})}\leq \left( \left\| \mathcal{I}(\Delta u_0)  \right\|_{L^{\infty}}
+\left\|  \mathcal{I}(\Delta u) \right\|_{L^{\infty}(0,t;L^{\infty})} +
\|\Delta \sigma\|_{L^2(0,t;L^2)}  \right).
\eeq
It remains to bound $\Delta\sigma$. Notice that we have
\[
\Delta \sigma= \mu \f{\p_y (\Delta u)}{\tau} + \mu (\Delta \vr)\p_y \widetilde{u}
-\left( Z^{\gamma_{+}}-(\widetilde{Z})^{\gamma_{+}} \right),
\]
and so, as for (\ref{bv7}) we obtain
\[
|\Delta \sigma| \leq C \Big( |\p_y (\Delta u)|+|\Delta \tau|(|\p_y \widetilde{u}|+1) +|\Delta R_0|+|\Delta Q_0|    \Big).
\]
It follows that
\eq{\label{bv10}
\|\Delta \sigma\|_{L^2(0,t;L^2)} \leq C \Big( \| \p_y (\Delta u)\|_{L^2(0,t;L^2)}+\|\Delta R_0\|_{L^{\infty}}+\|\Delta Q_0\|_{L^{\infty}}\\
+\int_0^t  (\|(\p_y \widetilde{u})(\cdot,s)\|_{L^2}+1) \|\Delta\tau(\cdot,s)\|_{L^{\infty}}  \ds  \Big).
}
Since  from (\ref{bm35}) we deduce that
$\int_0^T\|\p_y \widetilde{u}\|_{L^2}^2 \ds \leq C$, and therefore we can put together
(\ref{bv9})-(\ref{bv10}), and apply Gronwall's inequality to (\ref{bv8}) to deduce (\ref{bv4}). The proof of Lemma \ref{lem7} is thus finished. $\Box$

In order to use (\ref{bv4}) to conclude (\ref{bn19}), we need the estimates for $\Delta u$. In fact, standard energy estimate for parabolic equation \cite{EVS} gives
\begin{Lemma}\label{lem8}
For any $t\in (0,T]$, it holds that
\[
\|\Delta u\|_{L^{\infty}(0,t;L^2)}+\|\p_y(\Delta u)\|_{L^2(0,t;L^2)} \leq C\Big(\|\Delta u_0\|_{L^2} +
\|\Delta R_0\|_{L^{\infty}}+ \|\Delta Q_0\|_{L^{\infty}}
\]
\beq\label{bv12}
+ \|(\left\|(\p_y \widetilde{u})(\cdot,s)\|_{L^2}+1\right)\|\Delta\tau
(\cdot,s)\|_{L^{\infty}} \|_{L^2((0,t))}\Big).
\eeq
\end{Lemma}
Having this, (\ref{bn19}) follows by suitable combination of Lemmas \ref{lem7}-\ref{lem8}. For the sake of brevity, we omit the details and refer the reader  to \cite{LYS1} for similar steps. Clearly, (\ref{bn19}) implies the uniqueness of weak solutions and so the proof of Theorem \ref{lsz1} is complete. $\Box$

\section{Large time behavior of weak solution}\label{sion4}

In this section, we show the exponential decay of weak solution in $L^2$-norm. The classical methods to handle the large time behavior of the one-dimensional single-phase Navier-Stokes equations \cite{KZ,MY1,SV} are not readily applicable to our two-fluid model system. In \cite{Z1}, the author developed a new technique to treat one-dimensional viscous barotropic gas with nonmonotone pressure. Of great importance in \cite{Z1} is to obtain the uniform-in-time bounds of the density from above and below. It turns out that the idea can be adapted to our two-fluid model. As a matter of fact, it has already been successfully adapted before to the case of one-dimensional non-resistive magnetohydrodynamic equations \cite{LYS1}.

\subsection{Two-sided bounds for $R$ and $Q$}\label{sen41}

To begin with, we notice that the estimates in Lemma \ref{lem1} are uniform-in-time. Then we have the following lemma, which is essential for the proof of Theorem \ref{lsz2}. Throughout this section we use $C$ and $C_i$ to denote generic positive constants independent of time.
\begin{Lemma}\label{lem9}
Let $(R,Q,u)$ be the unique weak solution to (\ref{bn11})-(\ref{bn15}) ensured by Theorem \ref{lsz1}. Then
\beq\label{bh1}
C^{-1}\leq R(y,t),\,Q(y,t)\leq C,\,\, \text{ for a.e. } (y,t)\in [0,1]\times [0,\infty).
\eeq
\end{Lemma}
{\bf Proof.} From (\ref{bm44}) and the assumptions on the initial data (\ref{bn16})-(\ref{bn17}) one sees that verification of (\ref{bh1}) requires only to show the two-sided bounds for $\tau$. By adapting the arguments in \cite{Z1} (see also \cite{LYS1}), this follows from Lemma \ref{lem1} and the three items below.
\begin{itemize}
\item { $0 < C_1 \leq \int_0^1 \tau Z^{\gamma_{+}}\dy \leq C_2 <\infty$, }

\item { $Z^{\gamma_{+}}$ is sufficiently large if $\tau$ is sufficiently small,  }

\item { $Z^{\gamma_{+}}$ is sufficiently small if $\tau$ is sufficiently large.}
\end{itemize}
As a consequence, it remains to check that the three items above are satisfied. By the identity of pressure decomposition (\ref{bm8}) and (\ref{bm44}), it holds that
\eqh{
\int_0^1 \tau Z^{\gamma_{+}} \dy
&=\int_0^1 \left( \tau \alpha \left( \f{R}{\alpha}  \right)^{\gamma_{+}}
+\tau(1-\alpha)\left( \f{Q}{1-\alpha}  \right)^{\gamma_{-}} \right) \dy\\
&= \int_0^1 \Big(  (R_0\tau_0)^{\gamma_{+}}    (\alpha \tau)^{-\gamma_{+}+1} +
 (Q_0\tau_0)^{\gamma_{-}}    [(1-\alpha) \tau]^{-\gamma_{-}+1}  \Big) \dy\\
&\leq C_2,
}
where we have used the energy estimate (\ref{bm6}). Clearly, we conclude from the definition of $\alpha$, i.e., (\ref{bm4}), and Jensen's inequality that
\eqh{
\int_0^1 \tau Z^{\gamma_{+}} \dy
&=\int_0^1 \left( \tau \alpha \left( \f{R}{\alpha}  \right)^{\gamma_{+}}
+\tau(1-\alpha)\left( \f{Q}{1-\alpha}  \right)^{\gamma_{-}} \right) \dy\\
&\geq \int_0^1 \alpha^{-\gamma_{+}+1} \tau ^{-\gamma_{+}+1} (R_0\tau_0)^{\gamma_{+}}\dy\\
&\geq C \int_0^1 \tau ^{-\gamma_{+}+1}\dy\\
&\geq C \left(\int_0^1\tau \dy  \right)^{-\gamma_{+}+1}\\
&\geq C_1.
}
Suppose now that $\tau$ is small, i.e., $R+Q$ is large and we consider two possible cases. If $R$ is  large, then $Z^{\gamma_{+}}$ is also large due to $R\leq Z$.
If, on the other hand, $Q$ is large, then also $Z$ is large. Indeed, otherwise, we would arrive at a contradiction in the relation
\[
Q=\left(1-\f{R}{Z}\right) Z^{\gamma}.
\]
The third item is verified by using similar observation as above. We refer to \cite{Z1} and Lemma 5.3 in \cite{LYS1} for the remaining details.  $\Box$
\begin{Remark}
The key observations in Lemma \ref{lem9} are as follows. Firstly, the pressure term may be seen as a function with variables $y$ and $\tau$ by virtue of (\ref{bn9})-(\ref{bn10}), tending to infinity as $\tau$ goes to zero and tending to zero as $\tau$ goes to infinity. Secondly, the two internal pressures satisfy $\gamma$-laws. This leads to a positive lower bound of the integral $\int_0^1 \tau Z^{\gamma_{+}}\dy$; while the upper bound is obtained by the energy inequality. In this way, the arguments in \cite{LYS1,Z1} are naturally adapted.
\end{Remark}
\subsection{Exponential decay}\label{sen42}

In this subsection, we prove the exponential decay of weak solution in $L^2$-norm by adapting the ideas from \cite{Z1,LYS1}. It should be emphasized  that the structure of pressure function is crucial for a modification of these arguments to work.
\medskip

\noindent \emph{Step 1}. Let $(R_{\infty},Q_{\infty},u_{\infty})$ be the unique steady state for problem  (\ref{bn11})-(\ref{bn15}) given by (\ref{bn21}). Thanks to (\ref{bn21})$_3$, we rewrite the momentum equation (\ref{bn13}) as
\beq\label{bh2}
\p_t u +\p_y \left(Z^{\gamma_{+}}- Z_{\infty}^{\gamma_{+}  } \right) =\mu
\p_y \left( \f{\p_y u}{\tau} \right).
\eeq
As observed in the proof of Lemma \ref{lem2}, $Z$ can be seen as a function of $y$ and $\tau$, i.e., $Z=Z(y,\tau)$; whence $Z_{\infty}=Z(y,\tau_{\infty})$. Therefore, testing (\ref{bh2}) by $u$ and integrating by parts yields
\[
\f{1}{2}\f{d}{dt}\int_0^1 u^2 \dy+\int_0^1 \Big( Z^{\gamma_{+}}(y,\tau_{\infty})-Z^{\gamma_{+}}(y,\tau)   \Big)\p_y u\dy
+ \mu \int_0^1 \f{(\p_y u)^2}{\tau}\dy=0.
\]
Using the continuity equation (\ref{bn11}), one has
\eqh{
\int_0^1 \Big( Z^{\gamma_{+}}(y,\tau_{\infty})-Z^{\gamma_{+}}(y,\tau)   \Big)\p_y u\dy
&=\int_0^1 \Big( Z^{\gamma_{+}}(y,\tau_{\infty})-Z^{\gamma_{+}}(y,\tau)   \Big)\p_t \tau\dy\\
&=\f{d}{dt}\int_0^1 G(y,\tau,\tau_{\infty}) \dy,
}
where we denoted
\[
G(y,\tau,\tau_{\infty}):=\int_{\tau_{\infty}}^{\tau} \Big(Z^{\gamma_{+}}(y,\tau_{\infty})-
Z^{\gamma_{+}}(y,\xi)\Big)  d\xi.
\]
Thus we obtain
\beq\label{bh3}
\f{d}{dt}\int_0^1\left( \f{1}{2}u^2 +G(y,\tau,\tau_{\infty})   \right)  \dy+ \mu \int_0^1 \f{(\p_y u)^2}{\tau}\dy=0.
\eeq

\medskip

\noindent \emph{Step 2}. The key step in obtaining the exponential decay is to show that
\beq\label{bh4}
C^{-1}(\tau-\tau_{\infty})^2\leq G(y,\tau,\tau_{\infty}) \leq C(\tau-\tau_{\infty})^2.
\eeq
The main observation is as follows. By setting
\[
F(\tau):= -\int_{\tau_{\infty}}^{\tau} Z^{\gamma_{+}}(y,\xi) d\xi,
\]
$G(y,\tau,\tau_{\infty})$ is reformulated as
\beq\label{bh5}
G(y,\tau,\tau_{\infty})=F(\tau)-F(\tau_{\infty})-F'(\tau_{\infty})(\tau-\tau_{\infty}).
\eeq
Therefore, in order to deduce \eqref{bh4} it is enough to estimate the second derivative of $F(\tau)$.
To this purpose we use the expression for $\frac{\partial Z}{\partial \tau}$ from  (\ref{bm20}) to get
\eq{\label{bh6}
\f{\p Z^{\gamma_{+}}}{\p \tau}=\gamma_{+}Z^{\gamma_{+}-1}\f{\p Z}{\p \tau}
=-\gamma_{+}Z^{\gamma_{+}-1} \f{Q_0\tau_0\tau^{-2}+ R_0\tau_0\tau^{-2}Z^{\gamma-1}    }{\gamma Z^{\gamma-1}-R_0\tau_0\tau^{-1}(\gamma-1) Z^{\gamma-2}}.
}
As in (\ref{bm21}) we first observe that the denominator is strictly positive. Moreover, in spirit of Remark \ref{rk}, we infer from (\ref{bh1}) and the relation $Q=\left(1-\f{R}{Z}\right) Z^{\gamma}
$
that
\beq\label{bh7}
C^{-1}\leq Z(y,t)\leq C,\,\, \text{ for a.e. } (y,t)\in [0,1]\times [0,\infty),
\eeq
which together with lower and upper bound for $\tau$ implies boundedness of the numerator of \eqref{bh6}.
%As a consequence, making use of (\ref{bm21}), (\ref{bh1}) and (\ref{bh6})-(\ref{bh7}), we conclude from (\ref{bh5}) that (\ref{bh4}) is valid.

The remaining arguments follow largely the ones from \cite{Z1,LYS1}. We incorporate the detailed proof for the sake of completeness.

\medskip

\noindent \emph{Step 3}. Let $0<\ep <1$ and
\[
K(y,t):=\int_0^y (\tau(\xi,t)-\tau_{\infty}(\xi)) d\xi.
\]
Testing (\ref{bh2}) by $\ep K$ gives rise to
\eq{\label{bh8}
\f{d}{dt}\int_0^1 \varepsilon u K \dy-\varepsilon \int_0^1 \Big(Z^{\gamma_{+}}(y,\tau)  - Z^{\gamma_{+}}(y,\tau_{\infty})\Big)(\tau-\tau_{\infty})\dy\\
-\varepsilon \int_0^1 u^2 \dy +\varepsilon \int_0^1 \mu \f{\p_y u}{\tau} (\tau-\tau_{\infty})\dy=0.
}
From (\ref{bh3}) and (\ref{bh8}) we obtain
\[
\f{d}{dt}\int_0^1  \left(\f{1}{2}u^2 +G(y,\tau,\tau_{\infty})+\varepsilon u K\right)\dy
\]
\[
+\mu \int_0^1 \f{(\p_y u)^2}{\tau}\dy-\varepsilon \int_0^1 \Big(Z^{\gamma_{+}}(y,\tau)  - Z^{\gamma_{+}}(y,\tau_{\infty})\Big)(\tau-\tau_{\infty})\dy
\]
\beq\label{bh9}
=\varepsilon \int_0^1 u^2 \dy-\varepsilon \int_0^1 \mu \f{\p_y u}{\tau} (\tau-\tau_{\infty})\dy.
\eeq
It follows from (\ref{bh1}) and (\ref{bh6})-(\ref{bh7}) that
\beq\label{bh10}
\int_0^1 \Big(Z^{\gamma_{+}}(y,\tau)  - Z^{\gamma_{+}}(y,\tau_{\infty})\Big)(\tau-\tau_{\infty})\dy \geq C_1 \int_0^1 (\tau-\tau_{\infty})^2\dy.
\eeq
With the help of Cauchy-Schwarz's inequality and (\ref{bh1}), we find
\beq\label{bh11}
\left|\int_0^1 \mu \f{\p_y u}{\tau} (\tau-\tau_{\infty})\dy \right| \leq \f{C_2 }{2C_1}\mu \int_0^1 \f{(\p_y u)^2}{\tau}\dy +  \f{C_1 }{2} \int_0^1 (\tau-\tau_{\infty})^2 \dy;
\eeq
\beq\label{bh12}
  \int_0^1u^2\dy\leq C_3 \mu \int_0^1 \f{(\p_y u)^2}{\tau}\dy.
\eeq
Using (\ref{bh10})-(\ref{bh12}), (\ref{bh9}) implies
\[
\f{d}{dt}\int_0^1  \left(\f{1}{2}u^2 +G(y,\tau,\tau_{\infty})+\varepsilon u K\right)\dy
\]
\beq\label{bh13}
+\f{C_1 \varepsilon}{2} \int_0^1 (\tau-\tau_{\infty})^2 \dy+\left(1-\f{C_2 \varepsilon}{2C_1}-C_3\varepsilon \right)\int_0^1 \mu \f{(\p_y u)^2}{\tau} \dy \leq 0.
\eeq
\medskip

\noindent \emph{Step 4}. Due to the definition of $K$, it holds that
\beq\label{bh14}
\left| \int_0^1 \ep u K \dy  \right| \leq \f{\ep}{2}\int_0^1 u^2 \dy +\f{\ep}{2}\int_0^1 (\tau-\tau_{\infty})^2 \dy.
\eeq
Based on (\ref{bh14}), after choosing $\ep$ suitably small, we see
\[
C^{-1}(\|\tau-\tau_{\infty}\|_{L^2}^2+\|u\|_{L^2}^2)\leq
\int_0^1 \left(\f{1}{2}u^2 +G(y,\tau,\tau_{\infty})+\varepsilon u K\right) \dy\leq C (\|\tau-\tau_{\infty}\|_{L^2}^2+\|u\|_{L^2}^2),
\]
where we essentially used the property \eqref{bh4} from Step 2. Combining the above with (\ref{bh13}) leads to
\beq\label{bh15}
\|\tau-\tau_{\infty}\|_{L^2}+\|u\|_{L^2} \leq C \exp (-C t),
\eeq
for any $t\geq 0$.

\medskip

\noindent \emph{Step 5}. Finally, the exponential decay of $\|R-R_{\infty}\|_{L^2}$ and $\|Q-Q_{\infty}\|_{L^2}$ is a direct consequence of (\ref{bh15}) and the relations
\[
Q=Q_0\tau_0 \tau^{-1},\,\,R=R_0\tau_0 \tau^{-1};
\]
\[
Q_{\infty}=Q_0\tau_0 \tau_{\infty}^{-1},\,\,R_{\infty}=R_0\tau_0 \tau_{\infty}^{-1}.
\]
The proof of Theorem \ref{lsz2} is complete. $\Box$
\begin{Remark}
We observe that the exponential decay of $Z$ follows from that of $\tau$. Indeed,
\[
\|Z(y,\tau)-Z(y,\tau_{\infty})\|_{L^2}\leq \left\|\f{\p Z}{\p \tau}\right\|_{L^{\infty}}\|\tau-\tau_{\infty}\|_{L^2}\leq C \exp (-C t),
\]
in light of (\ref{bh6}), (\ref{bh7}) and (\ref{bh15}).
\end{Remark}

\begin{Remark}
The strategy adopted in this paper is strong enough  to show existence, stability and exponential decay of global weak solution to two-fluid models with more general form of pressure considered for example in \cite{NP}. In particular, the two-fluid model with pressure satisfying $\gamma$-laws \cite{LYS2} could be included.
\end{Remark}

\centerline{\bf Acknowledgement}
The research of Y. Sun is supported by NSF of China under grant number 11571167, 11771395, 11771206  and PAPD of Jiangsu Higher Education Institutions. The research of Y. Li is partially supported by  Postgraduate Research and Practice Innovation Program of Jiangsu Province under grant number KYCX 18-0028 and China Scholarship Council. E.Zatorska acknowledges the support from the Polish Ministry of Science and Higher Education ``Iuventus Plus"  under grant number 0888/IP3/2016/74.

%\abcdefg


\vskip0.5cm


\begin{thebibliography}{100}


\bibitem{AZ3} Amosov, A.A., Zlotnik, A.A.: {Uniqueness and stability of generalized solutions for a class of quasilinear systems of composite type equations.} { Math. Notes.} {\bf 55}, 555--567 (1994).

\bibitem{AKS} Antontsev, S.N., Kazhikhov, A.V., Monakhov, V.N.: {Boundary Value Problems in Mechanics of Nonhomogeneous Fluids}. North-Holland, Amsterdam, New York (1990).


\bibitem{HBV1} Beir\~{a}o da Veiga, H.: Long time behavior for one-dimensional motion of a general barotropic fluid. {Arch. Rational Mech. Anal.} {\bf 108}, 141--160 (1989).

\bibitem{BDGGH} Bresch, D., Desjardins, B., Ghidaglia,  J.M., Grenier, E., Hilliairet, M.:
Multifluid models including compressible fluids. Handbook of Mathematical
Analysis in Mechanics of Viscous Fluids, Eds. Y. Giga et A.
Novotn\'y, 1--52 (2017).

\bibitem{BrLi} Bresch, D.,  X. Huang, X.,  Li, J.: 
A Global Weak Solution to a One-Dimensional
Non-Conservative Viscous Compressible Two-Phase System. Comm. Math. Phys.
{\bf 309}, no. 3, 737--755 (2012).

\bibitem{BMZ} Bresch, D., Mucha, P.B., Zatorska, E: {Finite-energy solutions for compressible two-fluid stokes system.} {arXiv: 1709.03922, to appear in Arch. Ration. Mech. Anal. (2018).}

\bibitem{EVS} Evans, L.C.: {Partial Differential Equations.} 2nd ed., Graduate Studies in Mathematics, vol. 19, American Mathematical Society, Providence, RI, (2010).

\bibitem {EK1}Evje, S., Karlsen, K.H.: Global existence of weak solutions for a viscous two-phase model. { J. Diff. Equa. } {\bf 245}, 2660--2703 (2008).

\bibitem {EW1}Evje, S., Wen, H.: Weak solutions of a gas-liquid drift-flux model with general slip law for wellbore operators. { Discrete Contin. Dyn. Syst.} {\bf 33}, 4497--4530 (2013).

\bibitem {EW2}Evje, S., Wen, H.: Weak solutions of a two-phase Navier-Stokes model with a general slip law. { J. Funct. Anal.} {\bf 268}, 93--139 (2015).

\bibitem {EW3}Evje, S., Wen, H.: On the large time behavior of the compressible gas-liquid drift-flux model with slip. { Math. Models Methods Appl. Sci.} {\bf 25}, 2175--2215 (2015).

\bibitem {EW4}Evje, S., Wen, H., Yao, L.: Global solutions to a one-dimensional non-conservative two-phase model. { Discrete Contin. Dyn. Syst.} {\bf 36}, 1927-1-955 (2016).

\bibitem {EWZ1}Evje, S., Wen, H., Zhu, C.: On global solutions to the viscous liquid-gas model with unconstrained transition to single-phase flow. { Math. Models Methods Appl. Sci.} {\bf 27}, 323--346 (2017).



%\bibitem{FNP} Feireisl, E., Novotn\'{y}, A., Petzeltov\'{a}, H.: On the existence of globally defined weak solutions to the Navier-Stokes equations. { J. Math. Fluid Mech.} {\bf 3}, 358-392(2001)


\bibitem{IsHi} Ishii, M.,  Hibiki,  T.: Thermo-Fluid Dynamics of Two-Phase Flow.  Springer (2006).

\bibitem{KN1} Kanel, Y.: On a model system of equations for the one-dimensional motion of a gas. { Differ. Uravn.} {\bf 4}, 721--734 (1968).


\bibitem{KZ} Kazhikhov, A.V.: Stabilization of solutions of the initial-boundary value problem for barotropic viscous fluid equations. { Differ. Uravn.} {\bf 15}, 662--667 (1979).



\bibitem{LYS1} Li, Y., Sun, Y.: {Global weak solutions and long time behavior for 1D compressible MHD equations without resistivity.} {arXiv: 1710.08248}.


\bibitem{LYS2} Li, Y., Sun, Y.: {Asymptotic behavior of weak solutions to a compressible two-fluid model with large data.} {submitted}.


%\bibitem{LS} Lions, P.L.: {Mathematical Topics in Fluid Mechanics, Vol. 2, Compressible Models}. Clarendon Press, Oxford(1998)

\bibitem{MMM} Maltese, D., Mich\'{a}lek, M., Mucha, P.B., Novotn\'{y}, A., Pokorn\'{y}, M., Zatorska, E.: Existence of weak solutions for compressible Navier-Stokes equations with entropy transport. { J. Diff. Equa.} {\bf 261}, 4448--4485 (2016).


\bibitem{MY1} Matsumura, A., Yanagi, S.: Uniform boundedness of the solutions for a one-dimensional isentropic model system of compressible viscous gas. { Commun. Math. Phys.} {\bf 175}, 259-- 274 (1996).

\bibitem{NH} Nash, J.: {Le probl\'{e}me de Cauchy pour les \'{e}quations diff\'{e}rentielles d'un fluide g\'{e}n\'{e}ral,} {Bull. Soc. Math. France.} {\bf 90}, 487--497 (1962).

\bibitem{NP} Novotn\'{y}, A., Pokorn\'{y}, M.: {Weak solutions for some compressible multicomponent fluid models.} {arXiv: 1802.00798}.



\bibitem{SV} Stra\u{s}kraba, I., Valli, A.: Asympototic behavior of the density for one-dimensional Navier-Stokes equations. {Manuscripta Math.} {\bf 62}, 401--416 (1988).

%\bibitem{SZ1} Stra\u{s}kraba, I., Zlotnik, A.: On a decay rate for 1D-viscous compressible barotropic fluid equations. {J. Evol. Equ.} {\bf 2}, 69-96(2002)

\bibitem{VWY} Vasseur, A., Wen, H., Yu, C.: Global weak solution to the viscous two-fluid model with finite energy. { J. Math. Pure. Appl.} in press (2018).

\bibitem{WYZ} Wen, H., Yao, L., Zhu, C.:
Review on mathematical analysis of some two-phase flow models,
Acta Mathematica Scientia, {\bf 38}, 5,1617--1636 (2018).

%\bibitem{YZZ} Yao, L., Zhang, T., Zhu, C.: Existence and asymptotic behavior of global weak solutions to a 2D liquid-gas two-phase flow model. { SIAM J. Math. Anal.} {\bf  42}, 1874-1897(2010)


\bibitem{Z1} Zlotnik, A.A.: On equations for one-dimensional motion of a viscous barotropic gas in the presence of a body force. { Siberian Math. J.} {\bf 33}, 798--815 (1992).

\bibitem{ZA1} Zlotnik, A.A., Amosov, A.A.: Global generalized solutions of one-dimensional motion equations for a viscous barotropic gas. { Dokl. Akad. Nauk SSSR.} {\bf 299}, 1303--1307 (1988).




\end{thebibliography}
\end{document}